\documentclass[10pt,reqno,a4paper]{amsart}
\usepackage{amsmath,amsfonts,amssymb,amsthm,amsxtra}
\usepackage{diagrams}

\newtheorem{thm}{Theorem}[section]
\newtheorem{lem}{Lemma}[section]
\newtheorem{rem}{Remark}[section]
\newtheorem{prop}{Proposition}[section]
\newtheorem{exam}{Example}[section]
\newtheorem{defin}{Definition}[section]

\newcommand{\lL}{\mathbb{L}}

\newcommand{\T}{{\mathbb T}}

\newcommand{\R}{{\mathbb R}}

\begin{document}
\title{Extended Symmetries and Poisson Algebras Associated to Twisted Dirac Structures}
\author{Alexander Cardona}
\address{Mathematics Department
\\Universidad de Los Andes
\\A.A. 4976 Bogot\'a, Colombia.}
\email{acardona@uniandes.edu.co}
\date{April 17, 2012}

\maketitle

{\centerline{\textbf{To Steven Rosenberg, on his sixtieth
birthday.}}}

\begin{abstract}
In this paper we study the relationship between the extended
symmetries of exact Courant algebroids over a manifold $M$ defined
in \cite{BCG} and the Poisson algebras of admissible functions
associated to twisted Dirac structures when acted by Lie groups.
We show that the usual homomorphisms of Lie algebras between the
algebras of infinitesimal symmetries of the action, vector fields
on the manifold and the Poisson algebra of observables, appearing
in symplectic geometry,  generalize to natural maps of Leibniz
algebras induced both by the extended action and compatible moment
maps associated to it in the context of twisted Dirac structures.
\end{abstract}

\vspace{0.5cm}



\vspace{0.5cm}

\section{Introduction}
Let $\T M = TM \oplus T^* M$ denote the standard exact Courant
algebroid associated to a smooth manifold $M$, equipped with the
natural symmetric pairing
\begin{equation}\label{E:SPairing}
    \langle X \oplus \alpha, Y \oplus \beta \rangle
    = {1 \over 2} (  i_X \beta +  i_Y
    \alpha),
\end{equation}
where $X \oplus \alpha, Y \oplus \beta \in \Gamma (\T M)$, and the
twisted Dorfman bracket
\begin{equation}
  \label{E:TDB}
 [ X \oplus \alpha, Y \oplus \beta ]_H = [X, Y]\oplus \left(\mathcal{L}_X \beta  -
i_Y d \alpha - i_Y i_X H \right),
\end{equation}
where the twisting is given by the closed $3$-form $H$ on $M$. Let
$\lL_H < \T M$ be a \emph{Dirac structure} on $M$, i.e. a
sub-bundle of $\T M$ which is involutive under the bracket
(\ref{E:TDB}) and maximally isotropic with respect to $\langle
\cdot , \cdot \rangle$. The antisymmetrization of the bracket
(\ref{E:TDB}) gives rise to the twisted Courant bracket
\cite{C}\cite{SW}
\begin{equation}\label{E:TCB}
 [ X \oplus \alpha, Y \oplus \beta ]_H = [X, Y] \oplus \left( \mathcal{L}_X \beta
 -\mathcal{L}_Y \alpha - {1 \over 2}d (i_X \beta -
i_Y   \alpha) - i_Y i_X H \right),
\end{equation}
which, when evaluated on sections of $\lL_H$, coincides with the
twisted Dorfman bracket (\ref{E:TDB}).
 Twisted
Dirac structures appear naturally in Poisson geometry when, for
example, a reduction of a (twisted or non-twisted) Dirac structure
is performed \cite{BCG}. In quantum field theory and superstring
theory, the form $H$ has an interpretation as the Neveu-Schwarz
$3$-form \cite{Gra}. In \cite{C2012} it has been shown that,
associated to a twisted Dirac structure $\lL_H$, there is a
Poisson algebra of \emph{admissible functions} (the case of
non-twisted Dirac structures, studied by Courant and Weinstein in
\cite{C}\cite{CW}, is a particular case of this construction). In
general, a section $X \oplus \alpha \in \Gamma (\lL_H)$ is called
an admissible section, or \emph{admissible pair}, if \cite{C2012}
$$ d \alpha + i_X H =0,$$
and a smooth function $f$ on a manifold $M$ with a twisted Dirac
structure $\lL_H$ is called $H$-\emph{admissible} if there exists
a smooth vector field $X_f$ on $M$ such that $(X_f,df) \in
\Gamma(\lL_H)$ is an admissible pair, i.e. if $i_{X_f}H=0$. We
will denote by $C^\infty_{\lL_H}(M)$ the Poisson algebra of
$H$-admissible
functions on $M$ associated to $\lL_H$. \\
\\
In the case of Dirac structures associated to Poisson and
symplectic structures on $M$, which cannot be twisted, the set of
admissible functions is all of $C^\infty (M)$, but in general it
is not the case \cite{C2012}. If a function $f$ is $H$-admissible,
we will call a vector field $X_f$ such that $(X_f,df)$ is a
section of $\lL_H$ a \emph{Hamiltonian vector field} associated to
$f$. In \cite{C2012} it is shown that, in spite of the fact that
Hamiltonian vector fields are not unique in general, the bracket
\begin{equation}\label{E:PoissonBracket}
    \{ f, g \} = \mathcal{L}_{X_f}g
\end{equation}
defines a Poisson algebra structure on the space
$C^\infty_{\lL_H}(M)$ of $H$-admissible functions on $M$
(generalizing the classical result of \cite{C}). In this paper we
study the relation between the algebra of admissible functions in
the twisted case and the notion of moment map associated to
extended actions of Lie groups on exact Courant algebroids,
defined in \cite{BCG}. In particular we prove that extended
actions on Dirac structures, with compatible moment maps, induce
natural equivariant maps on the Lie algebra of vector fields and
the Poisson algebra of admissible functions associated to the
Dirac structure, giving rise to a relationship between Leibniz 
algebras and Poisson algebras of functions associated to Dirac
structures which generalize the known facts in symplectic and
Poisson geometry. \\
\\
The paper is organized as follows. In section \label{S:Admissible} we 
recall the notion of admissible pair for sections of Courant algebroids and 
Dirac structures, and the construction of the Poisson algebra associated to a 
twisted Dirac structure given in \cite{C2012}. In section \ref{S:Extended Actions} 
we recall the notions of Leibniz and Courant algebras, and we use them to extend some 
results in \cite{BCG} to the case of extended actions of Lie groups on exact Courant algebroids 
twisted by a closed $3$-form, together with the notion of moment map associated to such 
extended actions. 
In the last section we introduce the notion of Dirac actions and show, in theorem \ref{T:MuCourantAlgebra},
under which circumstances the usual morphisms of Lie algebras associated to 
Hamiltonian actions on symplectic manifolds can be recovered in this context, in
terms of morphisms of Leibniz algebras.

\vspace{0.3cm}

\section{Poisson Algebras Associated to Twisted Dirac Structures}\label{S:Admissible}
Let us consider, for $H \in \Omega^3(M)$ closed, a twisted Dirac
structure $\lL_H < \T M$  on $M$, i.e. a sub-bundle of $\T M$
which is involutive under the bracket (\ref{E:TDB}) and maximally
isotropic with respect to $\langle \cdot , \cdot \rangle$
\cite{SW}. As a first example consider the Dirac structure defined
by
\begin{equation}\label{E:SymplGraph}
    \lL_h =  \{(X, i_X h) \in \Gamma (\T M)\mid X \in \frak{X}(M)
    \},
\end{equation}
i.e. the graph in $\T M$ of a non-degenerate $2$-form $h$. It
follows from the definition of the twisted bracket (\ref{E:TDB})
that this Dirac structure is integrable if and only if $dh-H=0$,
so that $h$ cannot be closed in general (such a $h$ is called
 a $H$-closed $2$-form in \cite{SW}). Particular cases of Dirac
manifolds for which $H=0$ are Poisson and symplectic manifolds
(which correspond to graphs, in the generalized tangent bundle $\T
M$, of the corresponding Poisson bi-vector and symplectic form,
respectively). In these particular cases, the Poisson algebra
structure on $C^\infty (M)$ is defined by the action of
Hamiltonian vector fields on smooth functions given by
(\ref{E:PoissonBracket}). In general, even in the non-twisted
case, the Poisson algebra associated to a Dirac structure $\lL$ on
$M$ can be smaller than $C^\infty (M)$ since
(\ref{E:PoissonBracket}) defines a Poisson bracket only on
admissible functions associated to the Dirac structure, i.e. those
functions $f \in C^\infty (M)$ such that $(X_f, df) \in
\Gamma(\lL)$ for some $X_f \in \mathfrak{X}(M)$ (see e.g.\cite
{C}\cite{CW}). A further reduction is necessary in the case of
twisted Dirac structures \cite{C2012}.

\subsection{Admissible functions in the twisted case}
The notion of Dirac manifold $(M, \lL)$ as a sub-bundle of $\T M =
T M \oplus T^*M$ can be generalized to its higher analogues in
$\T^k M = T M \oplus \Lambda^k T^*M$, in the sense of \cite{Z},
where the integer $k \ge 0$ will be called the order of the Dirac
structure. Since $\Gamma (\T^k M) = \mathfrak{X} (M) \oplus
\Omega^k(M)$, the twisting in the corresponding Dorfman bracket
(\ref{E:TDB}) at order $k$ will be given by a closed $(k+2)$-form
on $M$. For $X \oplus \alpha, Y \oplus \beta \in \Gamma (\T^k M)$,
the twisted Dorfman bracket can be written as
\begin{equation}
 [ X \oplus \alpha, Y \oplus \beta ]_H =  \mathcal{L}_X Y
 + \mathcal{L}_X \beta  - i_Y (d \alpha + i_X H),
\end{equation}
where $H \in \Omega^{k+2}_{cl}(M)$, so that imposing $d \alpha +
i_X H=0$ is equivalent to impose a completely diagonal adjoint
action of $ X \oplus \alpha$ on $\Gamma (\T^k M)$:
\begin{equation}\label{E:Adjoint}
    \textsf{ad}_{ X \oplus \alpha} =\left(%
\begin{array}{cc}
  \mathcal{L}_X   & 0 \\
  0 & \mathcal{L}_X   \\
\end{array}%
\right).
\end{equation}
Actually, this is equivalent to consider the couple $(X , \alpha)$
as a geometric symmetry of the differential graded Lie algebra
associated to the dg-manifold $\textsf{Der}^\bullet(T[1]M \oplus
\R[k], Q_H)$, where the homological vector field is given by
\begin{equation}\label{E:QH}
    Q_H = d + H \partial_t,
\end{equation}
$d$ denotes the de Rham differential and $H \in
\Omega^{k+2}_{cl}(M)$ is the twisting. As a matter of fact, the
twisted Dorfman bracket (\ref{E:TDB}) is known to be the derived
bracket obtained from the complex of derivations
$\textsf{Der}^\bullet(T[1]M \oplus \R[k], Q_H)$ (see
\cite{Roy}\cite{S}\cite{U}). In \cite{C2012}, these facts are used
to motivate the following
\begin{defin}\label{D:AdmissiblePair}
Let $\lL_H$ in $\Gamma (\T^k M)$ be a $H$-twisted Dirac structure
of order $k$, where $H \in \Omega^{k+2}(M)$ is closed. A smooth
section $X \oplus \alpha \in \Gamma (\lL_H)$ is called an
admissible section, or \emph{admissible pair}, if
\begin{equation}\label{E:AdmissiblePair}
    d \alpha + i_X H =0.
\end{equation}
We will denote by $\Gamma_H (\T^k M)$ the space of admissible
pairs in $\Gamma (\T^k M)$.
\end{defin}
Notice that, when $k=0$, $\Gamma (\T^0 M) = \mathfrak{X}(M) \oplus
C^\infty (M)$ and interpreting the twisting $2$-form as a
symplectic form on $M$, condition (\ref{E:AdmissiblePair}) for a
section $(X,f)$ in $\Gamma (\T^0 M)$ is nothing but the usual
definition of the Hamiltonian vector field associated to $f$. When
$k=1$ equality (\ref{E:AdmissiblePair}) gives rise to a
\emph{exact derivation}  of the exact Courant algebroid $\T M$, in
the sense of \cite{BCG}. If $\lL_H$ is a twisted Dirac structure
in $\Gamma (\T^1 M)= \mathfrak{X}(M) \oplus \Omega^1 (M)$ it also
gives a criterium to identify a Poisson algebra of functions on
$M$ \cite{C2012}.
\begin{defin}\label{D:Admissible}
A smooth function $f$ on a manifold $M$ with a twisted Dirac
structure $\lL_H$ is called $H$-\emph{admissible} if there exists
a smooth vector field $X_f$ on $M$ such that $(X_f,df) \in
\Gamma(\lL_H)$ is an admissible pair, i.e. if $i_{X_f}H=0$. We
will denote by $C^\infty_{\lL_H}(M)$ the space of $H$-admissible
functions on $M$.
\end{defin}
If there is no twisting this definition of admissible function
coincides with the one of Courant in \cite{C}. On the other hand,
if the twisting is non-trivial, the set of $H$-admissible
functions may be smaller than the space of admissible functions in
the usual sense but, as shown in \cite{C2012}, it is a non-trivial
Poisson algebra.
\begin{thm}\label{T:Theorem}
Let $f,g$ be $H$-admissible functions on $M$ with respect to
the twisted Dirac structure $\lL_H$, where $H \in \Omega^3(M)$ is
closed. Then the product $fg$ and the bracket $\{f,g\}$ defined by
(\ref{E:PoissonBracket}) are $H$-admissible functions. Moreover,
such a bracket satisfies both Leibniz and Jacobi identities, and
then it defines a Poisson algebra structure on the space
$C^\infty_{\lL_H} (M)$.
\end{thm}
It is straightforward to see that restricting the twisted Dorfman
bracket (\ref{E:TDB}) to admissible pairs  $(X_f, df)$ and $(X_g,
dg)$ gives \cite{C2012}
\begin{equation}\label{E:AdmissBracket}
    [(X_f, df), (X_g, dg) ]_H =
([X_f,X_g] ,  d \{ f,g\}) ,
\end{equation}
generalizing the result already found in \cite{C} in the
non-twisted case. Moreover, since $$i_{[X_f,X_g] } H =
\mathcal{L}_{X_f } i_{X_g} H - i_{X_g}  \mathcal{L} _{X_f } H = -
i_{X_g} d i_{X_f } H + i_{X_g} i_{X_f } dH =0,$$ equation
(\ref{E:AdmissBracket}) implies that $\{ f, g\} $ is
$H$-admissible and
\begin{equation}\label{E:XBrackets}
   [X_f,X_g]  = X_{\{ f,g\}}.
\end{equation}
\begin{exam}\label{Ex:TSG}
Consider the Dirac structure defined in (\ref{E:SymplGraph}), i.e.
the graph in $\T M$ of a non-degenerate $2$-form $h$. It follows
from the definition of the twisted bracket (\ref{E:TDB}) that this
Dirac structure is integrable if and only if $dh-H=0$, and if the
functions $f, g, h \in C^\infty_{\lL_H} (M)$ are $H$-admissible,
$$  \{f,\{g,h \}\} +  \{g,\{h, f
\}\} +   \{h,\{f,g \}\} = H (X_f, X_g , X_h) =0, $$ so that Jacobi
identity holds. Notice that, being a graph of a non-degenerate
$2$-form, the twisted Dirac structures associates to any function
a Hamiltonian vector field $X_f$, but $f$ is a $H$-admissible
function through the condition $i_{X_f}H=0$ on $X_f$. Actually, a
pair $(X_f, df)$ in $\lL_h$ is $H$-admissible if and only if
$\mathcal{L}_{X_f} h = 0$. \\ \\ This Poisson algebra of functions
is not trivial in general: Let $(M, \omega)$ be  a symplectic
manifold and consider the $2$-form $h= \varphi \cdot \omega$,
where $\varphi \in C^\infty (M)$ has been chosen to make $h$
non-degenerate. Then the twisted Dirac structure
(\ref{E:SymplGraph}) is integrable with respect to the twisted
Courant bracket (\ref{E:TDB}) if and only if $H= dh = d \varphi
\wedge \omega$. For any $f\in C^\infty (M)$ there exists a vector
field $X_f \in \mathfrak{X}(M)$ such that $df = - i_{X_f} h$,
$(X_f, - df) \in \Gamma (\lL_H)$,  but  a smooth function $f$ is
$H$-a dmissible if and only if
$$ \mathcal{L}_{X_f} h =  \{f, \varphi \} \omega = 0.$$ Thus, in the cases in
which $\varphi$ is the Hamiltonian function for a dynamical system
with phase space $(M, \omega)$, an observable $f \in C^\infty (M)$
is $H$-admissible if and only if it is a constant of motion.
\end{exam}

This example suggests that condition (\ref{E:AdmissiblePair}) is a
symmetry condition when applied on functions. In section
\ref{S:Dirac Admissible and Symmetries} we will show that this is
the case since, actually, it combines both the requirement for an
action of a Lie group on $M$ to extend to an action on the exact
Courant algebroid $\T M$, and the requirements for this extended
action to have a compatible moment map associated to it.
\vspace{0.2cm}

\section{Extended actions and moment maps}\label{S:Extended Actions}
Given a smooth action of a Lie group $G$ on a manifold $M$, let us
denote by
\begin{equation}\label{E:InfAction}
\psi: \mathfrak{g} \to \mathfrak{X} (M)
\end{equation}
the associated infinitesimal action of $\mathfrak{g}$, the Lie
algebra of $G$, on $M$, which associates to each element $\xi$ in
$\mathfrak{g}$ the corresponding infinitesimal generator of the
action $X_{\xi} \in \mathfrak{X}(M)$. In this section we recall
the notion of Leibniz \cite{L} and Courant algebra, and the
definition of extension of the infinitesimal action to the exact
Courant algebroid $\T M$ \cite{BCG}. We show that the notion of
admissible pair given in (\ref{E:AdmissiblePair}) can be used to
characterize such extensions and also that
$\mathfrak{g}$-equivariant maps can be used to produce extensions
in the case of Courant algebroids twisted by an exact $3$-form.
The notion of moment map associated to an extended action
\cite{BCG} is also recalled.

\subsection{Extended actions of Lie groups}
In order to define an extension of the action $\psi: \mathfrak{g}
\to \mathfrak{X} (M)$ to sections of the exact Courant algebroid
$\T M= TM \oplus T^* M$,  the notion of Courant algebra was
introduced in \cite{BCG}. Recall that a \emph{Leibniz algebra}
$(\ell, [ \cdot , \cdot]_{\ell})$ is an algebra for which the
bilinear operation
$$ [ \cdot , \cdot]_{\ell} : {\ell} \times {\ell} \to {\ell}$$
is a derivation, i.e.
\begin{equation}\label{E:Leibniz}
[a  ,  [ b , c]_{\ell}]_{\ell}  = [[a , b]_{\ell} ,  c ]_{\ell}+
[b , [a , c]_{\ell}]_{\ell}
\end{equation}
for all $a  ,b ,c \in \ell $ \cite{L}. A morphism of  Leibniz
algebras is a homomorphism $f:  {\ell} \to {\ell'}$ such that
\begin{equation}\label{E:LeibnizMorphism}
f( [a  ,    b]_{\ell}) = [f(a)
, f(b)]_{\ell'}
\end{equation}
for all $a  ,b \in \ell $. A Leibniz algebra $(\ell, [ \cdot ,
\cdot]_{\ell})$ for which the bracket $ [ \cdot , \cdot]_{\ell} $
is antisymmetric is nothing but a Lie algebra. It follows that,
taking the quotient of a Leibniz algebra ${\ell} $ by the ideal
generated by the brackets of the form $ [ a , a]_{\ell} $, for all
$a \in {\ell} $, we obtain a Lie algebra $\mathfrak{g}_{\ell}$,
and the quotient map $f_{{\ell}}: {\ell} \to
{\mathfrak{g}}_{\ell}$ is a morphism of Leibniz algebras. A
natural way to build Leibniz algebras is considering
${\mathfrak{g}}$-modules and ${\mathfrak{g}}$-equivariant maps,
where ${\mathfrak{g}}$ denotes a Lie algebra. If
$\ell_{\mathfrak{g}}$ is a ${\mathfrak{g}}$-module and $\xi \cdot
\eta$ denotes the action of $\xi \in {\mathfrak{g}}$ on $\eta \in
\ell_{\mathfrak{g}}$, an application $\mu : \ell_{\mathfrak{g}}
\to {\mathfrak{g}}$ such that
$$ \mu (\xi \cdot \eta) = [\mu (\eta), \xi]_{\mathfrak{g}}$$
induces a Leibniz algebra structure on $\ell_{\mathfrak{g}}$ given by
\begin{equation}\label{Ex:Leibiz}
[\eta, \eta']_{\ell_{\mathfrak{g}}} = \mu(\eta') \cdot \eta,
\end{equation}
where $\eta, \eta' \in \ell_{\mathfrak{g}}$. Here again, the map
$\mu : \ell_{\mathfrak{g}} \to {\mathfrak{g}}$ is a morphism of
Leibniz algebras. Actually, every Leibniz algebra can be seen as
one of this type; this is the model of what is called a Courant
algebra in \cite{BCG} (see \cite{L} for more involved examples and
applications of Leibniz algebras).\\
\\
A \emph{Courant algebra} over a Lie algebra $\mathfrak{g}$ is a
Leibniz algebra $(\mathfrak{a}, [ \cdot , \cdot]_{\mathfrak{a}})$
with a morphism $\pi : \mathfrak{a} \to \mathfrak{g}$ of Leibinz
algebras, i.e.
\begin{equation}\label{E:CourantAlgebra}
 \pi([a,b]_{\mathfrak{a}}) = [\pi(a),
\pi(b)]_{\mathfrak{g}} ,
\end{equation}
for all $a,b,c \in \mathfrak{g}$.
 If $\pi$ is a surjective homomorphism and its kernel
$\mathfrak{h}= \ker \pi$ is abelian (with respect to $[ \cdot ,
\cdot]_{\mathfrak{a}}$) the Courant algebra is called
\emph{exact}, and in this case there is a natural
$\mathfrak{g}$-module structure on $\mathfrak{h}= \ker \pi$ given
by the adjoint action with respect to $[ \cdot ,
\cdot]_{\mathfrak{a}}$: If $\eta \in \mathfrak{h}$ and $a \in
\mathfrak{a}$ is such that $\pi (a) = \xi$ the map
\begin{equation}\label{E:gActionh}
\xi \cdot \eta = [a, \eta]_{\mathfrak{a}}
\end{equation} defines an action of
$\mathfrak{g}$ on $\mathfrak{h}$.
\begin{defin}\label{D:ExtendedAction} An extension of the action of a Lie group
$G$ on a manifold
$M$ to the Courant algebroid $\T M$ is an exact Courant algebra
$\mathfrak{a}$ over $\mathfrak{g}$, together with a Courant
algebra morphism $\rho: \mathfrak{a} \to \Gamma(\T M)$ such that
$\mathfrak{h}$ acts trivially and the induced action of
$\mathfrak{g}$ on $ \Gamma(\T M)$ integrates to a $G$-action on
$\T M$.
\end{defin}
Since the Courant algebroid $\T M$ we are working with is exact,
as noticed in \cite{BCG}, such an extension gives rise to a
commutative diagram
\begin{equation}\label{E:Diagram}
\renewcommand{\arraystretch}{1.5}
\begin{array}{cccccccccc}
&  0 & \xrightarrow{\text{ }} & \mathfrak{h} &\xrightarrow{\text{
}}& \mathfrak{a} &  \stackrel{\pi}{\longrightarrow}&   \mathfrak{g} & \xrightarrow{\text{ }}& 0\\
&  &  &\downarrow\rlap{$\scriptstyle \nu $} &  & \downarrow\rlap{$\scriptstyle \rho$} & &\downarrow\rlap{$\scriptstyle \psi$}& & \\
&  0 & \xrightarrow{\text{ }} & \Gamma (T^* M) &\xrightarrow{\text{ }}& \Gamma (\T M) &  \stackrel{\pi_{{}_{TM}}}{\longrightarrow}&   \Gamma (T M) & \xrightarrow{\text{ }}& 0\\
\end{array}
\end{equation}
in which the image of $\mathfrak{h} = \ker \pi$ under $\nu$ is
contained in $\Omega_{cl}^1(M)$ and, in order the action to
integrate to a $G$-action over $\T M$, we need a
$\mathfrak{g}$-invariant splitting of $\T M$. It turns out that
this condition is equivalent to ask the image of $\rho$ to be
given by admissible pairs in the sense of Definition
\ref{D:AdmissiblePair} (see also \cite{BCG}):
\begin{prop}\label{P:Extension} Let $G$ be a compact Lie group acting on a
smooth manifold $M$
and let $\pi: \mathfrak{a} \to \mathfrak{g}$ be an exact Courant
algebra with a morphism $\rho: \mathfrak{a} \to \Gamma(\T M)$ such
that $\nu(\mathfrak{h}) \subset \Omega^1_{cl} (M)$. Then $\rho$
extends to an action of the Courant algebra $\mathfrak{a}$ if and
only if $\rho (\mathfrak{a}) \subset \Gamma_{H}(\T M)$.
\end{prop}

\emph{Proof.}  Since $\nu (\mathfrak{h}) \subset
\Omega_{cl}^1(M)$, it follows from (\ref{E:TDB}) that
$\mathfrak{h}$ acts trivially, so we only have to verify that the
induced action of $\mathfrak{g}$ on $\T M$ integrates to an action
of $G$. It follows from (\ref{E:Adjoint}) that, if $\rho (a) \in
\Gamma_H (\T M)$ for all $a \in \mathfrak{a}$, the splitting in
$\T M = TM \oplus T^*M$ will be preserved and the action will
integrate to a $G$-action on $\T M$. Conversely, given an extended
action $\rho$, the usual averaging argument will give a
$\mathfrak{g}$-invariant splitting for $\T M$ $\Box$
\\

Among the Courant algebras over a Lie algebra $\mathfrak{g}$
induced by $\mathfrak{g}$-module structures, those induced by
semidirect products are particularly useful in order to define
extended actions. Consider a $\mathfrak{g}$-module $\mathfrak{h}$
with left action
$$ \cdot : \mathfrak{g} \times \mathfrak{h} \to \mathfrak{g}.$$
Restricting the adjoint action of $\mathfrak{g}$ on $\mathfrak{g}
\ltimes  \mathfrak{h}$ we have a Leibniz algebra, the
\emph{hemisemiditect product} of $\mathfrak{g}$ with
$\mathfrak{h}$, defined in \cite{KW}, which will be denoted  by
$(\mathfrak{a}_\mathfrak{g}^\mathfrak{h}, \cdot)$, with
multiplication given by
\begin{equation}\label{Hemi0}
(\xi,\eta)\cdot (\xi', \eta') = \left( [\xi,\xi'], \xi \cdot \eta'
\right) ,
\end{equation}
for all $(\xi,\eta),(\xi',\eta')  \in \mathfrak{g} \oplus
\mathfrak{h}$.

\begin{rem}\label{R:NewCourantAlgebra} Notice that an extended action
$\rho$ of the Courant algebra $(\mathfrak{a}, [ \cdot ,
\cdot]_{\mathfrak{a}})$ over $\mathfrak{g}$ on $\T M$ gives rise
naturally to a Courant algebra $\mathfrak{a}_M$ over
$(\mathfrak{X}(M), [\cdot , \cdot])$, given by $\mathfrak{a}_M=
\mathfrak{a}$ and
\begin{equation}\label{E:NewCourantAlgebra}
    \Pi_M : \mathfrak{a}_M \to \mathfrak{X}(M),
\end{equation}
where $\Pi_M (a) = \pi_{{}_{TM}}(\rho(a))$ for $a \in
\mathfrak{a}$.
\end{rem}

Recall that a map $ \upsilon : \mathfrak{h} \to \Omega^k(M)$
defined on a $\mathfrak{g}$-module $\mathfrak{h}$ is called
$\mathfrak{g}$-\emph{equivariant} if
\begin{equation}\label{E:Equivariance}
\upsilon (\xi \cdot \eta) = \mathcal{L}_{X_{\xi}}  \upsilon (\eta)
\end{equation}
for all $\xi, \eta \in \mathfrak{a}_\mathfrak{g}^\mathfrak{h}$.
The following proposition shows that equivariant maps give rise to
natural extensions of Lie algebra actions to hemisemidirect
product algebra actions on \emph{twisted} Courant algebroids.

\begin{prop}\label{P:TwistedExtension} Let $G$ be a compact Lie group
acting on a smooth manifold $M$
and let  $\mathfrak{h}$ be a $\mathfrak{g}$-module, where
$\mathfrak{g}$ denotes the Lie algebra of $G$. Given a
$\mathfrak{g}$-equivariant map $ \mu : \mathfrak{h} \to C^\infty
(M)$, the map   $\rho: \mathfrak{a}_\mathfrak{g}^\mathfrak{h} \to
\Gamma(\T M)$ given by
\begin{equation}\label{E:HdhExtension}
\rho (\xi, \eta) = (X_{\xi}, \alpha_{(\xi, \eta)}),
\end{equation}
where $X_{\xi}= \psi (\xi)$ and $ \alpha_{(\xi, \eta)}= d \mu
(\eta) + i_{X_{\xi}}h$, defines an extended action of the
hemisemidirect product  (\ref{Hemi0}) on the exact Courant
algebroid $\T M$, twisted by an exact $3$-form $H=dh$, if and only
if $\mathcal{L}_{X_{\xi}}h = 0$ for all $\xi \in \mathfrak{g}$.
\end{prop}

\emph{Proof.}  Consider $ (\xi, \eta),(\xi', \eta') \in
\mathfrak{a}_\mathfrak{g}^\mathfrak{h}$. Then, using Cartan
identities and (\ref{Hemi0}), we find that
\begin{eqnarray*}
 [\rho (\xi, \eta), \rho (\xi', \eta')]_{dh}
 &=& [X_{\xi}, X_{\xi'}] + \mathcal{L}_{X_{\xi}} ( d \mu (\eta') + i_{X_{\xi'}}h)
 - i_{X_{\xi'}} d (d \mu (\eta) + i_{X_{\xi}}h) - i_{X_{\xi'}}  i_{X_{\xi}} dh \\
 &=&X_{[{\xi}, {\xi'}]} + d \mu (\xi \cdot \eta') + i_{X_ {[{\xi}, {\xi'}] } }h   \\
  &=& \rho \left( (\xi, \eta)   \cdot (\xi', \xi')\right),
\end{eqnarray*}
since $\mu$ is $\mathfrak{g}$-equivariant. The result follows by
Proposition \ref{P:Extension}, since $d \alpha_{(\xi, \eta)} +
i_{X_{\xi}}H= \mathcal{L}_{X_{\xi}}h$ for all $\xi \in
\mathfrak{g}$
 $\Box$\\

Natural examples of extended actions on non-twisted Courant
algebroids include the actions commonly used in symplectic
geometry. If we consider $\mathfrak{h}= \mathfrak{g}$ in
(\ref{Hemi0}), and the adjoint action of $\mathfrak{g}$ on itself,
we obtain the exact Courant algebra
$\mathfrak{a}_\mathfrak{g}^\mathfrak{g}= \mathfrak{g} \oplus
\mathfrak{g}$ over $\mathfrak{g}$, with bracket
\begin{equation}\label{Hemi}
[(\xi,\eta), (\xi',\eta')]= \left( [\xi,\xi'], [\xi,\eta' ]
\right).
\end{equation}
\begin{exam}\label{Ex:Symplectic}
Let $M$ be  a smooth manifold and let $\omega$ be a closed
non-degenerate $2$-form on $M$. Consider the Courant algebroid $\T
M$ with $H=0$, and an infinitesimal action $ \psi : \mathfrak{g}
\to \mathfrak{X}(M)$ which integrates to a Lie group action on
$M$. The map $\rho: \mathfrak{a}_\mathfrak{g}^\mathfrak{g} \to
\Gamma (\T M)$ given by
\begin{equation}\label{E:HemiExtension}
\rho (\xi, \eta) = (X_{\xi}, \alpha_{\eta}),
\end{equation}
where $X_{\xi}= \psi (\xi)$ and $ \alpha_{\eta}=  i_{X_{\eta}}
\omega$, will give rise to an extended action whenever
$d\alpha_{\xi}=0$ for all $\xi \in \mathfrak{g}$. Indeed, as
follows from (\ref{Hemi}),
\begin{eqnarray*}
 [\rho (\xi, \eta), \rho (\xi', \eta')]
 &=&X_{[{\xi}, {\xi'}]} + i_{X_ {[{\xi}, {\eta'}] } } \omega   \\
  &=& \rho \left( [(\xi, \eta), (\xi', \eta')] \right),
\end{eqnarray*}
for all $\xi, \eta \in \mathfrak{g}$.
\end{exam}

\subsection{Moment Maps associated to Extended actions}
A moment map for an extended $\mathfrak{g}$-action $\rho:
\mathfrak{a} \to \Gamma(\T M)$ is a $\mathfrak{g}$-equivariant map
\begin{equation}\label{E:MomentMap}
    \mu : \mathfrak{h} \to C^\infty (M)
\end{equation}
such that $\nu = d \mu$, i.e. satisfying
\begin{equation}\label{E:Equivariance}
d \mu (\xi \cdot \eta) = \mathcal{L}_{\psi (\xi)} d\mu(\eta),
\end{equation}
where $ \psi :  \mathfrak{g} \to  \mathfrak{X}(M)$ denotes the
infinitesimal action of $G$ over $M$ and $\xi \cdot \eta$ is the
action given by (\ref{E:gActionh}). In \cite{BCG} the obstructions
to the existence of moment maps associated to extended actions
have been studied. It also has been shown that this definition of
moment map coincides with the usual one in symplectic geometry
when we consider the extended action given by
(\ref{E:HemiExtension}) on the Courant algebroid $\T M$, when
$H=0$ and $\omega$ is the symplectic form. This definition of
moment map is actually equivalent to ask the map $\mu \oplus \psi$
to induce an equivariant map $\rho_o: \mathfrak{a} \to \Gamma(\T^0
M) \cong \mathfrak{X}(M) \oplus C^{\infty}(M)$ such that the
following diagram commutes:

\begin{equation}\label{Diagram}
\begin{diagram}
&  &  0 &\rTo  & \mathfrak{h}&\rTo^{} &{\mathfrak{a}}&\rTo^{\pi} & \mathfrak{g} &\rTo  &0  \\
& &   & \ldTo(2,4)^{\;\;\;\;\;\;\;\;\; \mu} &  \dTo^{\nu}  & \ldTo(2,4)^{\;\;\;\;\;\;\;\;\; \rho_o}
&  \dTo^{\rho} & \ldTo(2,4)^{\;\;\;\;\;\;\;\;\; \psi} &  \dTo^{\psi} &   &   \\
& &0&\rTo &\Omega^1(M)&\rTo&\Gamma(\T^1 M)&\rTo& \mathfrak{X}(M) &\rTo& 0\\
&  & & \ruTo_d & & \ruTo && \ruTo_{id} & &   &&   \\
0 &\rTo &C^\infty(M)  &\rTo &\Gamma(\T^0 M)&\rTo& \mathfrak{X}(M)
&\rTo&0.
\end{diagram}
\end{equation}

\vspace{0.5cm} It is interesting to realize that, in this
approach, the moment map is no longer attached to the geometry,
i.e. to any particular Dirac structure in the exact Courant
algebroid $\T M$, but to the extended action itself. As a matter
of fact, example \ref{Ex:Symplectic} can be used to show that to
any equivariant map $\mu : \mathfrak{h} \to C^\infty (M)$, for a
$\mathfrak{g}$-module $\mathfrak{h}$, it is possible to associate
an extended action $\rho$ with moment map $\mu$ when $H=0$ (see
proposition 2.17 in \cite{BCG}). Proposition
\ref{P:TwistedExtension} before generalizes such result to the
twisted case when the twisting is \emph{exact}. In general, as we
will see in section \ref{S:Dirac Admissible and Symmetries}, the
existence of a moment map associated to an extended action amounts
to ``reduce" the space $C^\infty (M)$ in (\ref{Diagram}) to a
Poisson algebra of admissible functions with respect to the
twisting in $\T M$.

\begin{rem} In \cite{U} the equivariance of $\nu$ appears naturally
when a Hamiltonian action of the Lie group $G$ on the  dg-manifold
$\textsf{Der}^\bullet(T[1]M \oplus \R[k], Q_H)$  given in
(\ref{E:QH}) is defined as a map of Leibniz algebras
$$\mathfrak{g} \to \textsf{GDer}^\bullet(T[1]M \oplus \R[k], Q_H),$$
induced by the map of differential graded Lie algebras associated
to the infinitesimal action. Moreover, such maps of Leibniz
algebras are characterized in terms of invariant forms in the
Cartan model of equivariant cohomology
 (see \cite{U}, lemma 4.8).
\end{rem}

\vspace{0.3cm}

\section{Dirac structures, admissible functions and
symmetries}\label{S:Dirac Admissible and Symmetries} Let $\lL_H
\le \T M$ be a $H$-twisted Dirac structure on $M$, i.e. a
sub-bundle of the exact Courant algebroid $\T M$ which is
involutive under the bracket (\ref{E:TDB}) and maximally isotropic
with respect to the symmetric pairing (\ref{E:SPairing}). Consider
an extended action $\rho : \mathfrak{a} \to \Gamma(\T M)$
associated to an infinitesimal action $\psi : \mathfrak{g} \to
\mathfrak{X}(M)$ of a Lie group $G$ on $M$.
\begin{defin}\label{D:DiracAction}
The extended action $\rho$ will be called a \emph{Dirac action} on
$\lL_H$ if $\rho (a) \in \Gamma (\lL_H)$ for all $a \in
\mathfrak{a}$.
\end{defin}
Notice that the Dirac structure $\lL_H$ will be \emph{preserved}
by any Dirac action $\rho : \mathfrak{a} \to \Gamma(\T M)$ on it,
i.e. $[\rho (a) , \Gamma (\lL_H)]_H \subset \Gamma(\lL_H)$ for any
$a$ in the Courant algebra $\mathfrak{a}$. Dirac structures
preserved by extended actions give rise to reduced Dirac
structures \cite{BCG}. In this section we will show that, provided
the existence of moment maps, Dirac actions induce natural
equivariant maps between Courant algebras over the Lie algebra
$\mathfrak{g}$, giving rise to a relationship between Lie algebras
and Poisson algebras of functions associated to Dirac structures
which generalize the known facts in symplectic and Poisson
geometry.
\\

Let us first point out that, if $\rho (a) = (X_a, \alpha_a)$
defines a Dirac action, then the vector field $X_a$ should be a
symmetry of the twisting, i.e. $\mathcal{L}_{X_a}H=0$ for all $a
\in \mathfrak{a}$. As consequence of proposition \ref{P:Extension}
and remark \ref{R:NewCourantAlgebra} we have that both tangent and
cotangent components of a Dirac action $\rho$ are given by
equivariant maps:
\begin{lem}\label{L:Equivariance} Let $\rho (a) = (X_a, \alpha_a)$ denote a Dirac
action on a twisted Dirac structure $\lL_H$. Then, for any $a,b
\in \mathfrak{a}$,
\begin{equation}\label{E:DiracAction1}
    X_{[a,b]_\mathfrak{a}} = \mathcal{L}_{X_a} X_b
\end{equation}
and
\begin{equation}\label{E:DiracAction2}
    \alpha_{[a,b]_\mathfrak{a}} = \mathcal{L}_{X_a} \alpha_b.
\end{equation}
\end{lem}
Recall that in an exact Courant algebra $ \mathfrak{h}= \ker \pi
\to \mathfrak{a} \xrightarrow{\pi} \mathfrak{g}$ the action $\eta
\cdot \xi = [a, \xi]_{\mathfrak{a}}$ defines a
$\mathfrak{g}$-module structure on $\mathfrak{h}$, where $\eta \in
\mathfrak{h}$ and $a \in \mathfrak{a}$ is such that $\pi (a) =
\xi$. It follows then from (\ref{E:DiracAction1}) and
(\ref{E:DiracAction2}) that, in particular, the maps $X:
\mathfrak{h} \to \mathfrak{X}(M)$ and $\alpha : \mathfrak{h} \to
\Omega^{1}(M)$, defined by each component of the extended action,
are equivariant in the sense of (\ref{E:Equivariance}), i.e.
$X_{\xi \cdot \eta}= \mathcal{L}_{X_\xi} X_\eta $ and $\alpha_{\xi
\cdot \eta}= \mathcal{L}_{X_\xi} \alpha_\eta $ for all $\eta \in
\mathfrak{h}$ and $\xi \in \mathfrak{g}$. Actually
(\ref{E:DiracAction1}) was already observed in remark
\ref{R:NewCourantAlgebra}, and we will show that ---provided the
existence of a compatible moment map associated to the Dirac
action--- equation (\ref{E:DiracAction2}) induces a Courant
algebra structure over $(C^\infty_{\lL_H} (M), \{ \cdot,
\cdot\})$, the Poisson algebra given by theorem \ref{T:Theorem}.
\\
\\
An admissible function $f$ in the Poisson algebra
$C^\infty_{\lL_H} (M)$, associated to the $H$-twisted Dirac
structure $\lL_H \le \T M$, is a function for which there exists a
vector field $X_f$ such that $(X_f, df)$ is an admissible pair in
$\lL_H$, in the sense of definition \ref{D:Admissible}. If a Lie
group $G$ acts on $M$ by infinitesimal symmetries,  and such an
action extends to an action $\rho: \mathfrak{a} \to \Gamma (\T M)$
of a Courant algebra on the exact Courant algebroid $\T M$, we
have seen in proposition \ref{P:Extension} that $\rho (a)= (X_a,
\alpha_a)$ is an admissible pair for any $a \in \mathfrak{a}$,
i.e. $i_{X_a}H+ d \alpha_a =0$. For example, in the particular
case of the extended action $\rho: \mathfrak{g} \oplus
\mathfrak{g} \to \Gamma(\T M)$ given in example
\ref{Ex:Symplectic}, considering the Dirac structure $\lL_h$
associated to a non-degenerate $2$-form $h$ (non necessarily
closed, see example \ref{Ex:TSG}) given in (\ref{E:SymplGraph}),
with twisting $H=dh$, the condition on $\rho$ to be a Dirac action
implies that we have a ``diagonal" extended action: $\rho (a)=
(X_{\pi(a)} , i_{X_{\pi(a)}} h)$. Thus,
$$ \rho ([a, a']_{\mathfrak{g} \oplus
\mathfrak{g}})= [\rho (a), \rho (a')] =
 X_{[{\xi}, {\xi'}]} + i_{X_ {[{\xi}, {\xi'}] } } h,
$$
where $\pi(a)=\xi$ and $\pi(a')=\xi'$. Moreover, since $\rho (a)=
X_{\pi(a)} + i_{X_{\pi(a)}}h$ is an admissible pair, it follows
that the vector field $X_\xi = \psi(\xi)$ is ``locally
hamiltonian", i.e. $\mathcal{L}_{X_{\xi}} h = 0$ for any $\xi \in
\mathfrak{g}$. In this case a moment map $\mu: \mathfrak{g} \to
C^\infty(M)$ for such an extended action will give rise then to
\emph{admissible} functions $\mu_\xi \in C_{\lL_h}^\infty(M)$. If
$h$ is a closed form then for every smooth function $f$ on $M$
there exists a Hamiltonian vector field $X_f$ satisfying $i_{X_f}h
+ df =0$, so that $C_{\lL_h}^\infty(M)= C^\infty(M)$ and we have
the usual morphisms of Lie algebras in symplectic geometry:
\begin{equation}\label{D:SymplecticDiagram}
\begin{diagram}
&     &  &\mathfrak{g} & &     \\
&   & \ldTo^{\mu} &    & \rdTo^{\psi} &   & &    \\
&C^\infty(M)  & &\rTo^{\textsf{X}}& &\mathfrak{X}(M)
\end{diagram}
\end{equation}
associated to the infinitesimal action, where we have used
(\ref{E:XBrackets}). \\

In general, since a moment map associated to an extended action is
defined on the ``abelian part" of the Courant algebra $\pi:
\mathfrak{a} \to \mathfrak{g}$, this morphisms occurs very rarely.
Given a Dirac action $\rho: \mathfrak{a}_\mathfrak{g}^\mathfrak{g}
= \mathfrak{g} \oplus \mathfrak{g} \to \Gamma(\lL_H): a \mapsto
(X_a, \alpha_a)$ of the Courant algebra
$\mathfrak{a}_\mathfrak{g}^\mathfrak{g}$ defined in (\ref{Hemi})
on a twisted Dirac structure $\lL_H$ and a moment map $\mu:
\mathfrak{g} \to C^\infty(M)$ for the extended action, we will say
that such a moment map is \emph{compatible} with the action
whenever
\begin{equation}\label{E:CompatibleMu}
    \alpha_a = d \mu_{\bar{\pi}(a)},
\end{equation}
for all $a \in \mathfrak{a}$, where we denote by ${\bar{\pi}(a)}$
the pair $(0, \pi(a))$ in $\mathfrak{g} \oplus \mathfrak{g}$ in
order to distinguish it from $(0, {\pi}(a))$, for which
$\mu_{\pi(a)} = \mu ( {\pi}(a),0)=0$. In this case the morphisms
in diagram (\ref{D:SymplecticDiagram}) can be seen as particular
cases of the natural Leibniz algebra morphisms (Courant algebras)
induced both by the extended action and the moment map compatible
with it when $\mathfrak{a}=
\mathfrak{a}_\mathfrak{g}^\mathfrak{g}$.

\begin{thm}\label{T:MuCourantAlgebra} Let $\rho$ be
an extended action of the Courant algebra
$(\mathfrak{a}_\mathfrak{g}^\mathfrak{g}, [ \cdot ,
\cdot]_{\mathfrak{a}_\mathfrak{g}^\mathfrak{g}})$ on $\T M$, and
let $\lL_H$ be a twisted Dirac structure. If $\rho$ is a Dirac
action and there exists a moment map $\mu: \mathfrak{g} \to
C^\infty(M)$ compatible with it, then $\mu$ induces a Courant
algebra structure over $(C_{\lL_H}^\infty(M), \{\cdot , \cdot
\})$, given by:
\begin{equation}\label{E:MuCourantAlgebra}
    \Pi_\mu : \mathfrak{a}_\mu \to C_{\lL_H}^\infty(M),
\end{equation}
where $\mathfrak{a}_\mu= {\mathfrak{a}_\mathfrak{g}^\mathfrak{g}}$
and $\Pi_\mu (a) = \mu_{\bar{\pi} (a)}$ for $a \in
{\mathfrak{a}_\mathfrak{g}^\mathfrak{g}}$.\end{thm}

\textit{Proof.} Let $\rho(a)= (X_a, \alpha_a) \in \Gamma(\lL_H)$
denote the Dirac action of
${\mathfrak{a}_\mathfrak{g}^\mathfrak{g}}$ and let $\mu:
\mathfrak{g} \to C^\infty(M)$ be a moment map compatible with it.
Then $\mu (\eta) \in C_{\lL_H}^\infty(M)$ for all $\eta \in
\mathfrak{g}$ and
\begin{eqnarray*}
  \Pi_\mu ([a,b]_{{\mathfrak{a}_\mathfrak{g}^\mathfrak{g}}})  &=&   \mu_{[\pi(a), \pi(b)]}\\
    &=& \mathcal{L}_{X_\pi(a)} \mu_{\pi(b)},
\end{eqnarray*}
therefore, by (\ref{E:Equivariance}) in Lemma \ref{L:Equivariance}
and the definition (\ref{E:PoissonBracket}) of the Poisson bracket
in $C_{\lL_H}^\infty(M)$,
 \begin{equation}\label{E:PiMu}
\Pi_\mu ([a,b]_{{\mathfrak{a}_\mathfrak{g}^\mathfrak{g}}})=
\{\mu_{\pi(a)} , \mu_{\pi(b)} \},
\end{equation}
so that  $\Pi_\mu$ is a homomorphism of Leibniz algebras and
${\mathfrak{a}_\mathfrak{g}^\mathfrak{g}}$ is a Courant algebra on  $C_{\lL_H}^\infty(M)$
$\Box$
\\

Thus, when the hypothesis of theorem \ref{T:MuCourantAlgebra} are
fulfilled, we have a diagram of Leibniz algebra morphisms of the
form
\begin{equation}\label{CourantDiagram}
\begin{diagram}
&     &  &{{\mathfrak{a}_\mathfrak{g}^\mathfrak{g}}} & &     \\
&   & \ldTo(2,4)^{\Pi_\mu} &  \dTo^{\rho}  & \rdTo(2,4)^{\Pi_X} &   & &    \\
& &  &\Gamma(\T M)& & & & \\
& & \ldDotsto^{}& & \rdTo^{\psi} &&  & &   &   \\
&C_{\lL_H}^\infty(M)  & &\rTo^{\textsf{X}}& &\mathfrak{X}(M)
\end{diagram}
\end{equation}
attaching the Lie algebra $\mathfrak{X}(M)$ and the Poisson
algebra $C_{\lL_H}^\infty(M)$ to the Lie algebra of infinitesimal
symmetries $\mathfrak{g}$. The map $\rho_o$ in the lower row of
(\ref{Diagram}) ensembles the images of the Leibniz algebra maps
$\Pi_\mu$ and $\Pi_X$ as sections of $\T^0 M$. It is clear that
diagram (\ref{CourantDiagram}) becomes (\ref{D:SymplecticDiagram})
when $H=0$ and the Dirac structure is the graph
(\ref{E:SymplGraph}) of a closed non-degenerate $2$-form.

\begin{exam}\label{Ex:TSGSymmetry}
Consider the Dirac structure $\lL_h$ defined in
(\ref{E:SymplGraph}) as the graph in $\T M$ of the non-degenerate
$2$-form $h= \varphi \cdot \omega$ where $\omega$ denotes a
symplectic form on $M$. This Dirac structure is twisted by $H=d
\varphi \wedge \omega$, and $f \in C^\infty(M)$ is admissible if
and only if $\mathcal{L}_{X_f} h = \{f, \varphi \} \omega =0$.
Consider an action of a compact Lie group $G$ on $M$ such that
$\varphi$ is invariant, i.e. $\mathcal{L}_{X_\xi} \varphi =0$ for
all $\xi \in \mathfrak{g}$. Then the extended action
$$\rho (a) = (X_{\pi(a)}, i_{_{\bar{\pi}(a)}}h),$$
for $a \in {{\mathfrak{a}_\mathfrak{g}^\mathfrak{g}}}$ is a Dirac
action on $\lL_h$ with compatible moment map $\mu$. Since
$$\mathcal{L}_{X_\xi} h = 0$$
for $\xi \in \mathfrak{g}$ it follows that $\mu_\xi \in
C^\infty_{\lL_h} (M)$ for all $\xi \in \mathfrak{g}$, so that the
image of such a moment map is the set of ``constants of motion".
\end{exam}
 \vspace{0.9cm}

 {\bf Acknowledgements.}  The author is grateful to
Henrique Bursztyn, Michel Cahen, Simone Gutt, Yoshiaki Maeda and
Bernardo Uribe for many stimulating discussions on the geometry of
Poisson manifolds and Courant algebroids. This research has been
supported by the \emph{Vicerrector\'ia de Investigaciones} and the
\emph{Faculty of Sciences} of the Universidad de los Andes.

\end{document}